\documentclass{amsart}
\pdfoutput=1
\usepackage[T1]{fontenc}
\usepackage{amssymb,amsthm,amsmath,graphicx,pinlabel,wasysym,fullpage}
\usepackage[colorlinks]{hyperref}
\usepackage[style=default, margin=0pt, parskip=0pt, hangindent=0pt, indention=0pt, singlelinecheck=true, labelfont=up, justification=justified]{subcaption}

\newcommand{\R}{\mathbb{R}}

\newcommand{\MYhref}[3][blue]{\href{#2}{\color{#1}{#3}}}

\newcommand{\cH}{\mathcal{H}}

\newtheorem{theorem}{Theorem}[section]

\theoremstyle{definition}
\newtheorem{df}[theorem]{Definition}

\newtheorem{rmk}[theorem]{Remark}
\numberwithin{equation}{section}

\begin{document}

\renewcommand{\bf}{\bfseries}
\renewcommand{\sc}{\scshape}
\vspace{0.5in}

\title[Depicting a generalized shift move in crown diagrams]%
{Depicting a generalized shift move in crown diagrams}

\author{Jonathan D. Williams}
\address{Department of Mathematical Sciences, Binghamton University;
Binghamton, New York, 13902}
\email{jdw.math@gmail.com}


\begin{abstract} This paper gives a diagrammatic way to perform a generalized shift move on a crown diagram of a smooth 4-manifold. Applications include a simplified proof that if two crown diagrams are related by a generalized shift move, then they are slide-equivalent; a method for converting a genus $g>1$ Lefschetz fibration into a crown diagram; and the fact that the vanishing cycles of such a crown diagram are slide-equivalent to a standard inclusion of the Lefschetz vanishing cycles into a genus $g+1$ surface.
\end{abstract}

\maketitle

\section{\bf Introduction}
The study of the interplay between critical points of smooth mappings, the topology of their domains, and the diagrams that record them is a classical one which continues to yield interesting results, for example the relation between Lefschetz fibrations and trisections in \cite{CO}. This paper adds the corresponding result for a closely related object called a \emph{crown diagram}, by depicting how a rather local creation and cancellation homotopy involving the well-known critical points called \emph{indefinite cusps} can change the global fibration structure of a smooth 4-manifold as recorded in a crown diagram. The maps in this paper have a lot of names, such as \emph{excellent maps} \cite{S}, \emph{Morse 2-functions} \cite{GK}, and \emph{purely wrinkled fibrations} \cite{L}. Suffice it to say in this short introduction that they are maps from a smooth 4-manifold to a smooth 2-manifold whose singular sets consist of indefinite cusps and indefinite folds. The paper then uses this knowledge to give an algorithm to turn a Lefshetz fibration into a crown map while keeping track of the vanishing cycles.

In \cite{W1} the author introduced the notion of crown diagrams. A crown diagram is a pair $(\Sigma,\Gamma)$ where $\Sigma$ is a smooth, connected, closed, oriented 2-manifold without boundary, and $\Gamma=\{c_i:i=1,2,\ldots,\ell\}$ is a collection of simple closed curves in $\Sigma$ satisfying certain conditions which are mostly not relevant to this paper, but which are spelled out in \cite[Remark 1.6]{W2} following ideas from \cite[Section 2.2]{H}. In that same remark, it is explained why the oriented diffeomorphism class of a crown diagram specifies a smooth, connected, closed, oriented 4-manifold $M$ without boundary. Crown diagrams come from \emph{crown maps} $M\to S^2$, and most of the rest of \cite{W2} is devoted to proving that if two crown maps are homotopic, then their corresponding crown diagrams are related by certain moves that can change $\Sigma$ or $\Gamma$. Diagrammatic interpretations of these moves exist for all but one: the shift move. It is the aim of this paper to supply it. 

A previous treatment of the cusp merge homotopy appeared in \cite{BH}; that paper describes all versions of the shift move and of merge homotopies in the most generality possible, resulting in descriptions in the language of mapping class groups. This paper has a more explicit approach, giving a diagrammatic description of a cusp merge and a generalized shift chosen from the infinitude of possible merges and shifts, obtained diagrammatically, and observes in Remark \ref{rmk2.1} that the other choices are related to the one that is given by multislide, whose effect on diagrams is easier to diagrammatically understand than shift.

To explain the term \emph{slide-equivalent}, a \emph{slide} on crown diagrams is visually the same as a handleslide move on Heegaard diagrams of 3-manifolds:

\begin{df}Given two closed curves $a$ and $b$ in $\Sigma$, one \emph{slides} $a$ over $b$ by replacing $a$ with a connected sum $a^s=a\#b$.\end{df} 

This differs from Heegaard handleslides in a few respects; for example $a^s$ might not be a simple closed curve, depending on $a$, $b$, and the path in $\Sigma$ between them used to form the connected sum. On the other hand, the slides in this paper occur in sequences which result in closed curves which are homotopic to simple closed curves.

\begin{df}This paper says $(\Sigma,\Gamma)$ is \emph{slide-equivalent} to $(\Sigma,\Gamma')$ when there is a sequence of slides and isotopies of individual circles applied to $(\Sigma,\Gamma)$ that results in $(\Sigma,\Gamma')$.\end{df}

The fact that slide equivalence is an equivalence relation on $\ell$-element sequences of closed curves in a fixed surface for each $\ell$ follows from the fact that the inverse of a slide is a slide.

\begin{figure}
	\begin{subfigure}{0.30\textwidth}
		\centering
		\labellist
		\small\hair 2pt
		\pinlabel $f$ at 24 55
		\pinlabel $f$ at 60 55
		\pinlabel $f$ at 84 55
		\endlabellist
		\includegraphics[width=\textwidth]{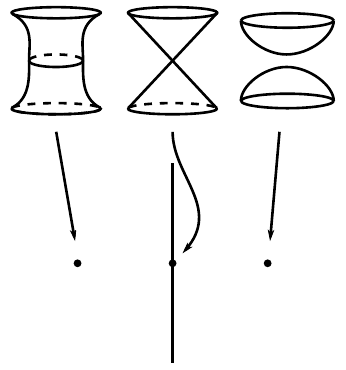}
		\caption{The vertical line is the image of a fold arc.}
		\label{fold}
	\end{subfigure}
	\hspace{0.04\textwidth}
	\begin{subfigure}{0.30\textwidth}
		\centering
		\includegraphics[width=\linewidth]{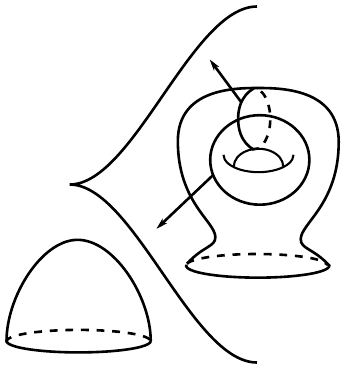}
		\caption{The image of a cusp point between two fold arcs.}
		\label{cusp}
	\end{subfigure}
	\hspace{0.04\textwidth}
	\begin{subfigure}{0.30\textwidth}
		\centering
		\labellist
		\small\hair 2pt
		\pinlabel $p$ at 85 23
		\pinlabel $p_1$ at 77 66
		\pinlabel $p_2$ at 50 85
		\pinlabel $p_3$ at 21 66
		\pinlabel $p_4$ at 20 37
		\pinlabel $q$ at 47 46
		\pinlabel $p_{\ell-1}$ at 55 20
		\pinlabel $\gamma$ at 66 46
		\endlabellist
		\includegraphics[width=\linewidth]{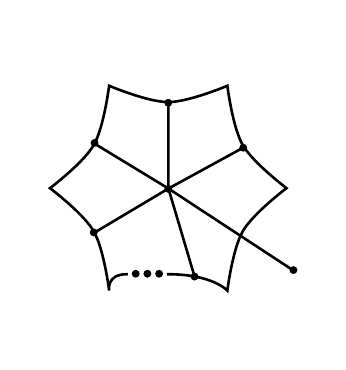}
		\caption{The critical image of a crown map, decorated with a graph $T$.}
		\label{crown}
	\end{subfigure}
	\caption{Figures \ref{fold} and \ref{cusp} represent the behavior of the fibers over points near the critical image of a crown map $f$. As a point moves from left to right in Figure \ref{fold}, a circle in the fiber collapses to a point, then the two cones separate, increasing the Euler characteristic of the fiber by 2. Figure \ref{cusp} illustrates how the circles associated to the two folds that meet at a cusp must intersect transversely at one point. In that figure, each fiber is drawn in the region that contains its image.}\label{crit}
\end{figure}

Crown maps and the more general \emph{indefinite Morse 2-functions} have appeared in many places in the literature. Polynomial models and diagrammatic interpretations of their critical points and those of their deformations appear in \cite{B1,GK,L,W1} and many other papers, though such information is not required for this paper. Figures \ref{fold} and \ref{cusp} are included to establish the visual vernacular of the paper and remind the reader of what the fibration near a fold arc or a cusp point looks like.

\section{\bf Surgered diagrams of Morse 2-functions which are injective on their critical sets}\subsection{Introduction to surgered diagrams}\label{intro2diags} In Figure \ref{crown}, the critical image of the crown map forms a cusped circle, and one may use a connection $\cH$ and the pictured reference paths from $q$ to $p_i$, $i=1,2,\ldots,\ell$ to mark the vanishing cycles $\Gamma=\{c_1,\ldots,c_\ell\}$ as simple closed curves in the fiber over $q$. This is where the crown diagram $(\Sigma,\Gamma)$ comes from.

Instead of using the fiber over $q$, one could depict the vanishing cycles in the fiber $F_p$ over $p$ using the reference paths given by following $\gamma$ from $p$ to $q$, then following each of the $\ell-1$ reference paths coming from $q$ to the points $p_1,\ldots,p_{\ell-1}$. There are two points in $F_p$ that flow along $\gamma$ to the critical point in the fiber where $\gamma$ crosses the $\ell^{th}$ fold arc, and the corresponding vanishing cycle $c_\ell$ appears as a neighborhood of that pair of points in $F_p$, which this paper will depict as a pair of shaded disks, essentially corresponding to the two disks at the right of Figure \ref{fold}. The vanishing cycles that intersected $c_\ell$ in the previous paragraph now appear as arcs in $F_p$ connecting the two disks that form this neighborhood. It will be useful later to use two disks in $F_p$ instead of two points in order to more easily keep track of vanishing cycles as reference paths change (see for example Figure \ref{surgcusp2}), so that $(\Sigma,\Gamma)$ can be recovered from the resulting collection of arcs at every moment of a shift deformation. This is a \emph{surgered} diagram as in \cite[Section 6.1]{BH} and \cite[Remark 1.7]{W2}. The circles bounding the pair of disks can be thought of as the two ends of a tube that, when attached to $F_p$, form the fiber $F_q$ over $q$. This pair of circles can also be thought of as two copies of the vanishing cycle $c_\ell$.

In the case that the Morse 2-function has more than one critical circle, each of which is the boundary of a disk of regular values on its higher-genus side like in Figure \ref{cuspmergefig1}, a surgered diagram is simply the union of vanishing cycles and shaded regions from a surgered diagram for each circle.

This paper assumes familiarity with the results of \cite[Section 4]{BH}, which is a treatment of facts such as the uniqueness of a vanishing cycle as measured using reference paths that differ by postcomposing with an ambient isotopy in $S^2$, or using connections that differ by homotopy (Theorem 4.6), and the interpretation of a surgered diagram whose reference path is as in the bottom center of Figure \ref{surgcusp2} (Lemma 4.7). This paper adds slightly to the latter story in the case that the reference path is a member of a family of reference paths (the addition is slight: the pair of surgered diagrams at either side of the bottom of Figure \ref{surgcusp2}).

\subsection{How a surgered diagram changes when its reference path moves across a cusp}\label{surgcuspsection} 
\begin{figure}
	\begin{subfigure}[t]{0.4\textwidth}
		\centering
		\labellist
		\small\hair 2pt
		\pinlabel $\gamma_0$ at 80 90
		\pinlabel $p$ at 112 107
		\pinlabel $q$ at 50 107
		\pinlabel $1$ at 48 122
		\pinlabel $\ell$ at 48 92
		\pinlabel $c_\ell$ at 23 43
		\pinlabel $c_\ell$ at 70 43
		\pinlabel $c_1$ at 40 48
		\pinlabel $b$ at 47 60
		\pinlabel $g$ at 5 47
		\pinlabel $g$ at 88 47
		\endlabellist
		\includegraphics[width=\textwidth]{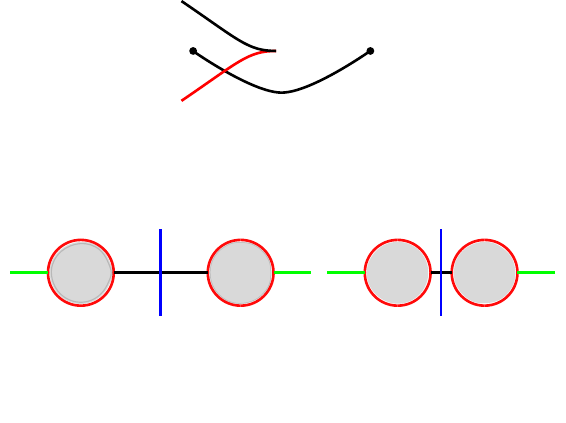}
		\caption{Top: The critical image near a cusp of a crown map, decorated with reference points $p,q$ and a reference path $\gamma_0$.\\Bottom left: part of the corresponding surgered diagram.\\Bottom right: The reference path has moved slightly toward the cusp, causing $c_1$ to shorten.
		}
		\label{surgcusp1}
	\end{subfigure}
	\hspace{0.1\textwidth}
	\begin{subfigure}[t]{0.4\textwidth}
		\centering
		\labellist
		\small\hair 2pt
		\pinlabel $\gamma_{1/2}$ at 93 100
		\pinlabel $c_1,c_\ell$ at 102.5 42
		\endlabellist
		\includegraphics[width=\linewidth]{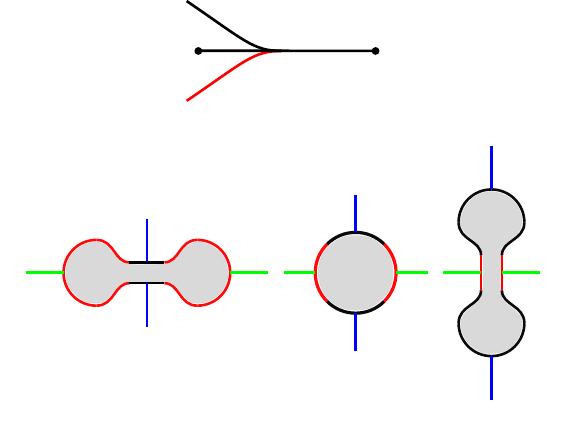}
		\caption{Top: The reference path crosses the cusp.\\Bottom: Three diagrams show how the surgered diagram evolves by ambient isotopy in the lower-genus fiber $F_p$ before the reference path passes the cusp.
		}
		\label{surgcusp2}
	\end{subfigure}
	
	\vspace{.25cm}
	
	\begin{subfigure}{0.4\textwidth}
		\centering
		\labellist
		\small\hair 2pt
		\pinlabel $\gamma_1$ at 83 124
		\pinlabel $c_1$ at 96 43
		\pinlabel $c_1$ at 141.5 43
		\pinlabel $g$ at 167 43
		\pinlabel $b$ at 118 75
		\pinlabel $b$ at 118 10
		\pinlabel $c_\ell$ at 112 48
		\endlabellist
		\includegraphics[width=\linewidth]{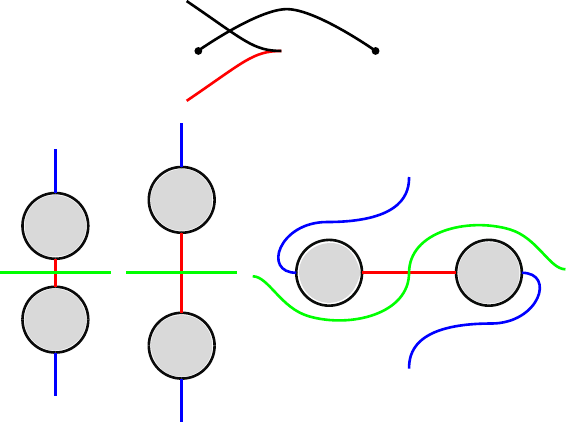}
		\caption{Top: The reference path has passed the cusp.\\Bottom: The surgered diagram undergoes a further ambient isotopy so that the union of disks and $c_\ell$ coincides with the earlier arrangement.
		}
		\label{surgcusp3}
	\end{subfigure}
	\caption{This figure consists of three pictures of the image of the critical set near a cusp in $S^2$ above the relevant parts of eight corresponding surgered diagrams. As the family of reference paths $\gamma_t$ from $p$ to $q$ crosses the cusp, the surgered diagrams evolve so that the same crown diagram is recoverable at every point in time. Neighborhoods of points are shaded, with colors corresponding to those of fold arcs. The arcs $g$ and $b$ are parts of hypothetical vanishing cycles.}\label{surgcusp}
\end{figure}

On the top of Figure \ref{surgcusp1} there is a cusp where two fold arcs meet in the critical image of a crown map, two reference points $p$ and $q$, and a reference path $\gamma_0$ connecting them. The reference path $\gamma_0$ is the initial member of a family of reference paths $\gamma_t$, $t\in[0,1]$. At the bottom left of Figure \ref{surgcusp1}, with labels $g,c_\ell,c_1$ and $b$, appears a decorated disk in the fiber above $p$. This disk is a neighborhood of the pair of vanishing cycles $c_1$ and $c_\ell$ in the surgered diagram. The decorations include a shaded neighborhood labeled $c_\ell$ of the pair of points that flow to the lower fold arc $\ell$. The vanishing cycle $c_1$ of the upper fold arc 1 appears as a line segment connecting the two shaded regions $c_\ell$ because the local model for a cusp requires the two vanishing cycles $c_1$ and $c_\ell$ to intersect at one transverse point in the fiber. The figure includes a few other segments to represent where other vanishing cycles $g$ and $b$ may intersect $c_1$ and $c_\ell$. Note $g$ and $b$ may be contained in the same vanishing cycle, or a vanishing cycle might intersect $c_\ell$ or $c_1$ multiple times, in which case that vanishing cycle would correspond to multiple parallel copies of $g$ or $b$. Directly to the right of this picture is another copy of it in which $c_1$ has shortened. This corresponds to a movement of the reference path so that it crosses the critical image at a point which is closer to the cusp, but still crosses the lower fold arc $\ell$.

As the crossing between the reference path and the fold arc $\ell$ begins to move upward toward the cusp, $c_1$ shortens in the diagram below until, at Figure \ref{surgcusp2}, both $c_1$ and $c_\ell$ correspond to the single point of $F_p$ that flows above $\gamma_{1/2}$ from $F_p$ to the cusp point in $M$. For a  surgered diagram with a pair of disks like in Figure \ref{surgcusp1}, the boundary of the pair of disks corresponds to the boundary of a neighborhood of $c_\ell$ in the higher genus fiber $F_{q}$. For a reference path like $\gamma_{1/2}$ that crosses a cusp like in Figure \ref{surgcusp2}, the boundary of the single disk corresponds to the boundary of a neighborhood of $c_1\cup c_\ell$, which is a genus-1 connect summand of $F_{q}$. For this reason, the disk should be better understood to have four sides, with alternating sides corresponding to $c_1$ and $c_\ell$ as suggested by the colors in the figure. The validity of this interpretation can be seen by the fact that the crown diagram remains recoverable up to diffeomorphism from the surgered diagrams throughout the evolution depicted in Figure \ref{surgcusp}, at every instant using the same correspondence between crown diagrams and surgered diagrams. With this understood, the change to Figure \ref{surgcusp3} is an application of the same principle in reverse. The rightmost picture of the fiber in Figure \ref{surgcusp3} is included to illustrate one of the many ways the shaded regions can be moved to coincide with their original position in case a cusp merge is the next step, as is the case in Section \ref{shiftmove}.

\subsection{How a cusp merge effects a surgered diagram}\label{cuspmergesection} 
\begin{figure}
	\begin{subfigure}{0.45\textwidth}
		\centering
		\labellist
		\small\hair 2pt
		\pinlabel $p$ at 85 57
		\pinlabel $\rightarrow$ at 175 30
		\pinlabel $q_1$ at 29 30
		\pinlabel $q_2$ at 140 30
		\pinlabel $\ell_1$ at 66 40
		\pinlabel $k_1$ at 66 20
		\pinlabel $\ell_2$ at 100 40
		\pinlabel $k_2$ at 100 20
		\endlabellist
		\includegraphics[width=\linewidth]{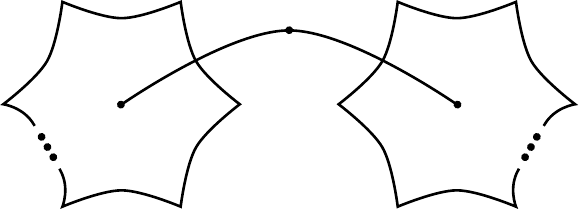}
		\caption{\ }
		\label{cuspmergefig1}
	\end{subfigure}\qquad
	\begin{subfigure}{0.45\textwidth}
		\centering
		\includegraphics[width=\linewidth]{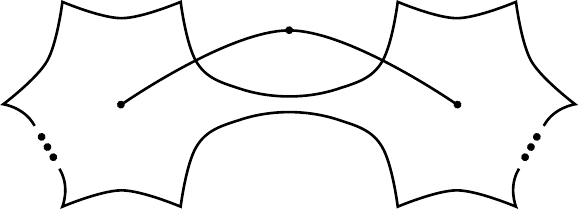}
		\caption{\ }
		\label{cuspmergefig2}
	\end{subfigure}
	\caption{Cusp merge.}\label{cuspmergefig}
\end{figure}
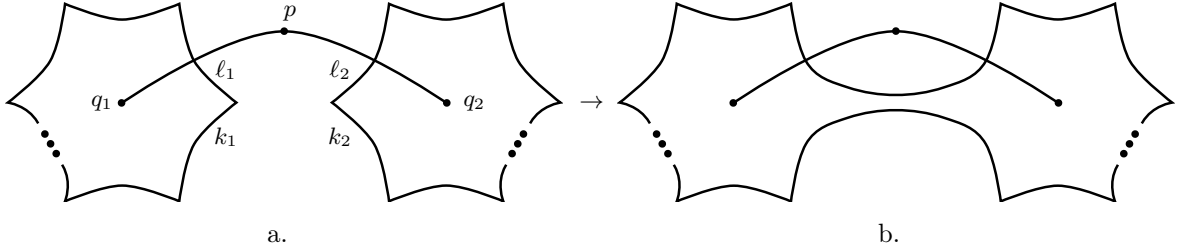

Figure \ref{cuspmergefig} depicts how the image of the critical set changes during a cusp merge. In this case, the merge involves the two cusps between $\ell_i$ and $k_i$, $i=1,2$. This homotopy was first discussed in the context of 4-manifold topology as \emph{Move 2} in \cite[Section 3]{L} as a homotopy of maps $\R^4\to\R^2$. This paper uses Lekili's result that such a homotopy exists, given a Morse 2-function whose critical image is as in Figure \ref{cuspmergefig1}, and finds a crown diagram for Figure \ref{cuspmergefig2}. 

\begin{rmk}\label{rmk2.1}As discussed in \cite{BH} the cusp merge is performed according to a choice of \emph{framed joining cuve} in $M$ from one cusp to another, and this choice can greatly affect the resulting crown diagram. On the other hand, the fibrations resulting from two choices of cusp merge path are homotopic: The homotopy is simply the reverse of one of the cusp merge homotopies (such a homotopy is called a fold merge), followed by the other cusp merge. As defined in \cite[Section 2.3]{W2}, this homotopy corresponds to a multislide move on crown diagrams, so the crown diagrams that result from two choices of generalized shift are related by a multislide. This concludes the remark.\end{rmk}

Any surgered diagam for Figure \ref{cuspmergefig1} is the union of the \emph{stationary set} and \emph{nonstationary set} in $F_p$, slightly generalizing terminology from \cite{W2}. Each of these sets is the collection of disks and arcs in a surgered diagram for one of the critical circles in the picture. That is, each set has a pair of shaded disks corresponding to a fold arc $\ell_i$ along with the vanishing cycles for the rest of the fold arcs going around $q_i$. Recall that the vanishing cycle $k_i$ must appear as a single embedded arc from one disk corresponding to $\ell_i$ to the other disk corresponding to $\ell_i$. To explain the terms \emph{stationary} and \emph{nonstationary}, suppose the nonstationary set is chosen to correspond to the circle around $q_1$. There is an ambient isotopy of the nonstationary set in $F_p$ sending the pair of disks corresponding to $\ell_1$ to coincide with those of $\ell_2$, and also sending the vanishing cycle of $k_1$ to coincide with that of $k_2$. This isotopy preserves the condition that the union of stationary and nonstationary sets in $F_p$ forms a surgered diagram for the fibration, and performing the cusp merge homotopy also preserves the validity of the vanishing cycles, though the two coinciding depictions of $c_{\ell_i}$ and the coinciding pair of vanishing cycles $c_{k_i}$ now are each double depictions of single vanishing cycles which might as well be called $c_\ell$ and $c_k$, respectively. The result is a surgered diagram for Figure \ref{cuspmergefig2}, from which a crown diagram can be immediately recovered. As in Section \ref{surgcuspsection}, the validity of this paragraph simply comes from the fact that crown diagrams for each circle are recoverable from the surgered diagrams at every moment, whether stationary, nonstationary, or combined as in the last step. 

To paraphrase the previous paragraph: Given a surgered diagram corresponding to a pair of reference paths like in Figure \ref{cuspmergefig1}, isotope the nonstationary set so that the relevant pair (a vanishing cycle $c_{k_i}$ connecting the pair of disks $\ell_i$) coincides with that of the stationary set, then delete one of the now-redundant pairs and draw the corresponding crown diagram. The ordering of the vanishing cycles after the cusp merge comes immediately from the ordering of the corresponding fold arcs of the crown map.

Given a surgered diagram derived using reference paths as in Figure \ref{cuspmergefig1}, a cusp merge between a different pair of cusps can be achieved as in the previous two paragraphs after using Section \ref{surgcuspsection} to change which fold arcs are crossed by the reference paths.

\section{Application: Generalized shift move}\label{shiftmove} This section gives a diagrammatic way to apply the shift move of \cite[Section 2.4]{W1}. It leads to a much simpler and more self-contained proof of a statement that first appeared in \cite[Remark 4.5]{W3}, which is as follows:
\begin{theorem}\label{thm1}If the crown diagram $(\Sigma,\Gamma')$ is obtained by applying a generalized shift move to the crown diagram $(\Sigma,\Gamma)$, then $(\Sigma,\Gamma)$ is slide-equivalent to $(\Sigma,\Gamma')$.\end{theorem}
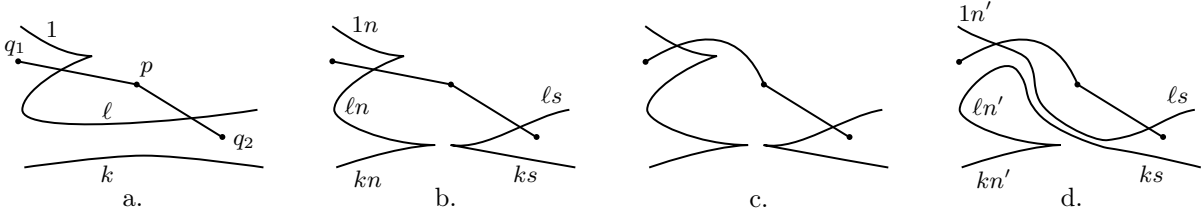
\begin{figure}
	\begin{subfigure}{0.2\textwidth}
		\centering
		\labellist
		\small\hair 2pt
		\pinlabel $1$ at 14 53
		\pinlabel $\ell$ at 35 22
		\pinlabel $k$ at 35 -2
		\pinlabel $p$ at 50 38
		\pinlabel $q_1$ at 0 46
		\pinlabel $q_2$ at 87 10
		\endlabellist
		\includegraphics[width=\linewidth]{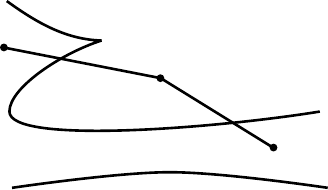}
		\caption{\ }
		\label{shiftdef1}
	\end{subfigure}\hspace{.045\textwidth}
	\begin{subfigure}{0.2\textwidth}
		\centering
		\labellist
		\small\hair 2pt
		\pinlabel $1n$ at 14 55
		\pinlabel $\ell n$ at 11 24
		\pinlabel $kn$ at 15 -4
		\pinlabel $\ell s$ at 85 28
		\pinlabel $ks$ at 75 -2
		\endlabellist
		\includegraphics[width=\linewidth]{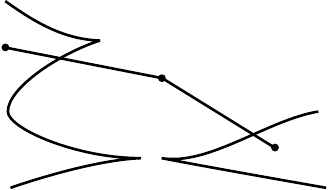}
		\caption{\ }
		\label{shiftdef2}
	\end{subfigure}\hspace{.045\textwidth}
	\begin{subfigure}{0.2\textwidth}
		\centering\includegraphics[width=\linewidth]{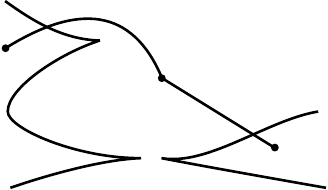}
		\caption{\ }
		\label{shiftdef3}
	\end{subfigure}\hspace{.045\textwidth}
	\begin{subfigure}{0.2\textwidth}
		\centering
		\labellist
		\small\hair 2pt
		\pinlabel $1n'$ at 8 60
		\pinlabel $\ell n'$ at 13 24
		\pinlabel $kn'$ at 15 -4
		\pinlabel $\ell s$ at 85 28
		\pinlabel $ks$ at 75 -2
		\endlabellist
		\includegraphics[width=\linewidth]{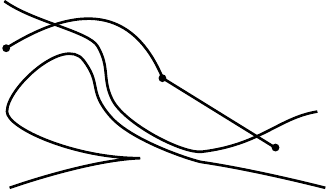}
		\caption{\ }
		\label{shiftdef4}
	\end{subfigure}
	\caption{The shift deformation consists of a fold merge followed by a cusp merge. The letters $n$ and $s$ designate fold arcs corresponding to nonstationary and stationary vanishing cycles, respectively.}\label{shiftdef}
\end{figure}
\begin{proof}For the crown map whose critical image is partly depicted in Figure \ref{shiftdef1}, the surgered diagrams corresponding to the two reference paths from $p$ to $q_1$ and from $p$ to $q_2$ can be assumed to be exactly the same. Having performed the fold merge, the two diagrams now coexist in the surgered diagram for Figure \ref{shiftdef2} as discussed in the last paragraph of Section \ref{intro2diags}. At this point there are two sets of vanishing cycles, and the ones whose reference paths contain $q_1$ will be the nonstationary set. For this reason, following the labels in Figure \ref{shiftdef2} counter-clockwise, the indices of the nonstationary set will have the suffix $n$:
\[c_{1n},\ c_{2n},\ \ldots,\ c_{kn},\  c_{\ell n}\] 
while the stationary set will be written $c_{ks},\ c_{(k+1)s},\ \ldots,\ c_{\ell s}$.
	
To obtain a surgered diagram for Figure \ref{shiftdef3}, apply the modification of Figure \ref{surgcusp} to the nonstationary set, resulting in a list of vanishing cycles 
\[c_{1n}',\ c_{2n}',\ \ldots,\ c_{kn}',\ c_{\ell n}'\]
corresponding to the same respective fold arcs as the un-primed list.
	
Performing the cusp merge following Section \ref{cuspmergesection}, the surgered diagram undergoes the following identifications:
\[c_{1n}'=c_{\ell s},\hspace{1cm}c_{\ell n}'=c_{ks}.\] The first equation involves two copies of each vanishing cycle and the pair of disks they bound, while the second involves one vanishing cycle on each side of the equation. Thus the cusp merge results in the following sequence of vanishing cycles, which can be read off from the labels of Figure \ref{shiftdef4}: 
\begin{equation}\label{pvc1}
	c_{1n}',\ c_{2n}',\ \ldots,\ c_{kn}',\ c_{\ell n}',\ c_{(k+1)s},\ c_{(k+2)s},\ \ldots,\ c_{(\ell-1)s}.\end{equation}
The goal is to prove this sequence is slide-equivalent to $c_1,c_2,\ldots,c_\ell$. To this end, observe that Figure \ref{surgcusp3} depicts how the nonstationary vanishing cycles change between Figures \ref{shiftdef2} and \ref{shiftdef3}. This implies that moving the reference path across the cusp causes $c_{1n}$ and $c_{\ell n}$ to switch places in the surgered diagram; that is,
\[c_{1n}=c_{\ell n}',\hspace{1cm}c_{\ell n}=c_{1n}'.\]Also recall that the non-primed vanishing cycles are equal to the original vanishing cycles; that is, $c_{i\ast}=c_i$ for any $i\in\{1,\ldots,\ell\}$ and for any $\ast\in\{n,s\}$ such that $c_{i\ast}$ was defined above. Accordingly, the sequence (\ref{pvc1}) can be rewritten
\begin{equation}\label{pvc2}c_\ell,\ c_{2n}',\ \ldots,\ c_{kn}',\ c_1,\ c_{(k+1)s},\ c_{(k+2)s},\ \ldots,\ c_{(\ell-1)s}.\end{equation}

If an element of $\{c_{2n},\ldots,c_{kn}\}$ does not intersect $c_{1n}$ or $c_{\ell n}$, then it is unchanged by the modification in Figure \ref{surgcusp}, so its prime and its suffix may be removed. Since the stationary set was defined as the set which is unchanged, each suffix $s$ may also be removed. Since the ordering of vanishing cycles is cyclic, $c_\ell$ can be moved to the end of the list:
\begin{equation}\label{pvc3}c_{2n}',\ \ldots,\ c_{kn}',\ c_1,\ c_{k+1},\ c_{k+2},\ \ldots,\ c_{\ell-1},c_\ell.\end{equation}

If an element of $\{c_{2n},\ldots,c_{kn}\}$ intersects $c_{1n}$ or $c_{\ell n}$, then the intersection can be depicted as that between $b$ or $g$ and $c_{1n}$ or $c_{\ell n}$ in Figure \ref{surgcusp1}, and the change in $b$ and $g$ from the very first diagram in Figure \ref{surgcusp1} to the very last diagram in Figure \ref{surgcusp3} can be realized by sliding over $c_{1n}$ and over $c_{\ell n}$ one time each. For this reason, the sequence (\ref{pvc3}) is slide-equivalent to
\begin{equation}\label{pvc4}c_2,\ c_3,\ \ldots,\ c_k,\ c_1,\ c_{k+1},\ c_{k+2},\ \ldots,\ c_\ell.\end{equation}
As observed in \cite[Remark 4.5]{W3}, it is straightforward to move $c_1$ from its place in the sequence \ref{pvc4} to the beginning of the sequence by repeatedly transposing it with the vanishing cycle that immediately precedes it using a sequence of slides.\end{proof}

\section{Application: Crown diagrams from Lefschetz fibrations}\label{leftocrown}
\begin{figure}
\begin{subfigure}[t]{0.45\textwidth}
	\centering
	\labellist
	\small\hair 2pt
	\pinlabel $p$ at 70 110
	\pinlabel $q_i$ at 112 110
	\pinlabel $l_i$ at 103 50
	\endlabellist
	\includegraphics[width=\linewidth]{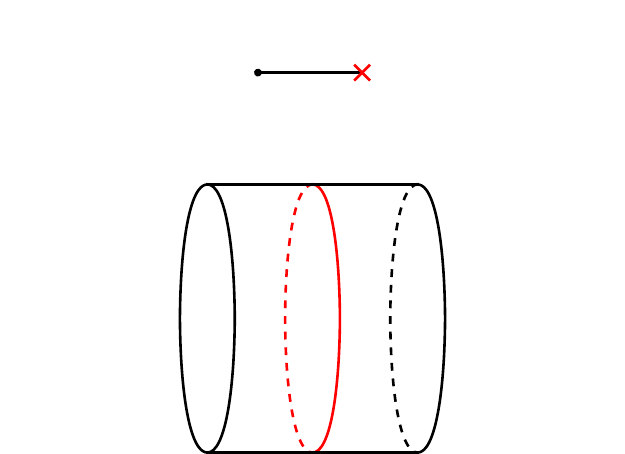}
	\caption{Top: A reference path from the regular value $p$ to the the Lefschetz critical value $q_i$.\\Bottom: A neighborhood of the Lefchetz vanishing cycle $l_i$ in the preimage $F_p$ of $p$.}
	\label{wrinklefig1}
\end{subfigure}\hspace{.045\textwidth}
\begin{subfigure}[t]{0.45\textwidth}
	\centering
	\labellist
	\small\hair 2pt
	\pinlabel $r_i$ at 105 120
	\pinlabel $g_i$ at 105 99
	\pinlabel $r_i$ at 101.5 5
	\pinlabel $g_i$ at 105 37.5
	\pinlabel $b^1_i$ at 98 53
	\pinlabel $b^2_i$ at 98 21.5
	\endlabellist
	\includegraphics[width=\linewidth]{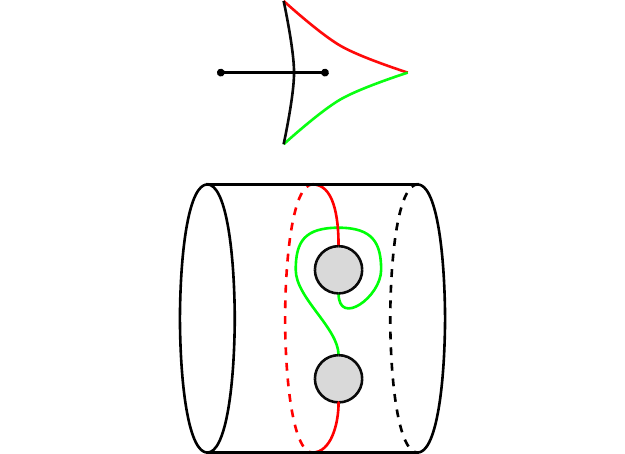}
	\caption{Top: The image of a wrinkled Lefschetz critical point, with the same reference path ($q_i$ is now a regular value).\\Bottom: The vanishing cycles of the three fold arcs. The labels $r_i,b^j_i,g_i$ correspond to $a,b,d$ in \cite[Figure 6]{L}, respectively.}
	\label{wrinklefig2}
\end{subfigure}\vspace{.3cm}
\begin{subfigure}[t]{0.7\textwidth}
	\centering
	\labellist
	\small\hair 2pt
	\pinlabel $r_1$ at 110 79
	\pinlabel $r_n$ at 132 79
	\pinlabel $g_i$ at 103 56
	\pinlabel $b^1_i$ at 121 56.5
	\pinlabel $b^2_i$ at 121 24.5
	\pinlabel $r_1$ at 110 0
	\pinlabel $r_n$ at 132 0
	\endlabellist
	\includegraphics[width=\linewidth]{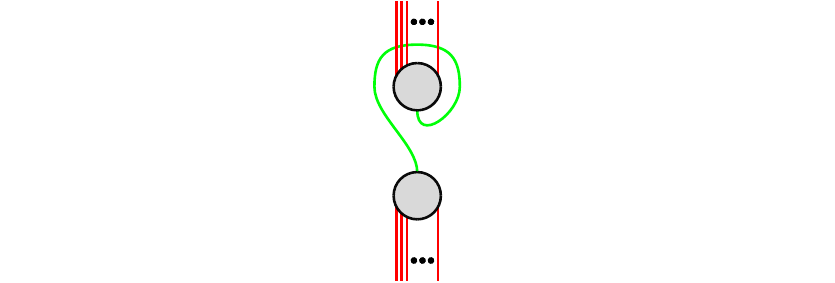}
	\caption{A region of the surgered diagram for the first group of triangles in Figure \ref{mergesfig}. Outside of this region, the diagram coincides with $r_1,\ldots,r_n$ in $F_p$. Having wrinkled all Lefschetz critical points, the curves $g_i,\ b^1_i$ and $b^2_i$ coincide, respectively, up to isotopy. }
	\label{wrinklefig3}
\end{subfigure}\hspace{.045\textwidth}
\caption{A way to obtain a surgered diagram for a Morse 2-function which comes from wrinkling the critical points of a Lefschetz fibration $M\to S^2$. Figure \ref{wrinklefig3} depicts a pseudocornoation if $g_i$ and $b_i^j$ are omitted for all $i\neq1$.}\label{wrinklefig}
\end{figure}

A smooth Lefschetz fibration $M\to S^2$ is commonly recorded as a sequence of simple closed curves $l_1,\ldots,l_n$ in the preimage $F_p$ of a regular value $p$ of the fibration. The curves $l_i$ are the vanishing cycles of the Lefschetz critical points of the fibration as measured using reference paths $\gamma_1,\ \gamma_2,\ \ldots, \gamma_n$, where $\gamma_i$ goes from $p$ to the Lefschetz critical value $q_i$ \cite[Chapter 8]{GS}.
In \cite[Move 4]{L} there appeared a way to convert a Lefschetz fibration into a Morse 2-function by perturbing the map near its critical points, thereby converting each Lefschetz critical point as in Figure \ref{wrinklefig1} into a critical circle with three cusps as in Figure \ref{wrinklefig2}. Each triangle has its own trio of vanishing cycles and pair of disks as in the figure, and up to isotopy of the diagram the disks $b_i^1$ coincide, the disks $b_i^2$ coincide, and the vanishing cycles $b_i,$ that form their boundaries coincide, for $i\in\{1,\ldots,n\}$.
\begin{df}\label{corodef}For a Lefschetz fibration $M\to S^2$ given by a sequence of Lefschetz vanishing cycles $l_1,\ldots,l_n$, a \emph{coronation} of $l_1,\ldots,l_n$ is the corresponding list of vanishing cycles, starting with $b_1$, of a crown map obtained by wrinkling the critical points of the fibration and merging the resulting triangles. A \emph{coronation of a Lefschetz fibration} is a coronation of any such sequence of its Lefschetz vanishing cycles. A \emph{pseudocoronation} is the sequence of simple closed curves $b_1,r_1,g_1,r_2,\ldots,r_n$ in the surgered diagram for the Morse 2-function obtained from wrinkling the Lefschetz critical points as in Figure \ref{wrinklefig}.\end{df}
\begin{rmk}Usually, when discussing the vanishing cycles of a Lefschetz fibration, the Lefschetz critical points are shown to be mapped along a circle, with the corresponding vanishing cycles listed counter-clockwise. It is ultimately a consequence of the order of vanishing cycles in Figure \ref{wrinklefig2} that the sequence $b_1,r_1,g_1,r_2,\ldots,r_n$ in Definition \ref{corodef} is a \emph{clockwise} sequence, and it is slide-equivalent to the corresponding clockwise sequence of vanishing cycles, starting with $b_1$, of the crown map.\end{rmk}
\begin{theorem}\label{thm2}Any coronation of a sequence of Lefschetz vanishing cycles $l_1,\ldots,l_n$ in a closed surface of genus $g>1$ is slide-equivalent to any pseudocoronation of $l_1,\ldots,l_n$.\end{theorem}

Though this theorem is about Lefschetz fibrations over the sphere, its proof also applies word-for-word to fibrations over a disk; however this is beyond the scope of the paper. The following will be useful in the proof.

\begin{df}Given two sequences of circles $a_1,\ldots,a_n$ and $b_1,\ldots,b_n$ in a surface are slide-equivalent, $a_i$ is said to be slide-equivalent to $b_i$ for $i\in\{1,\ldots,n\}$. Given a sequence of circles $a_1,\ldots,a_n$ in a surface, the notation $a_i^s$ is used to denote any closed curve which is slide-equivalent to $a_i$.\end{df}

\begin{proof}Let the $i^{th}$ triangle be the one whose vanishing cycles are $b_i,g_i,r_i$, and let the \emph{main circle} begin as the first triangle, into which all other triangles will be merged one by one as in Figure \ref{mergesfig}.
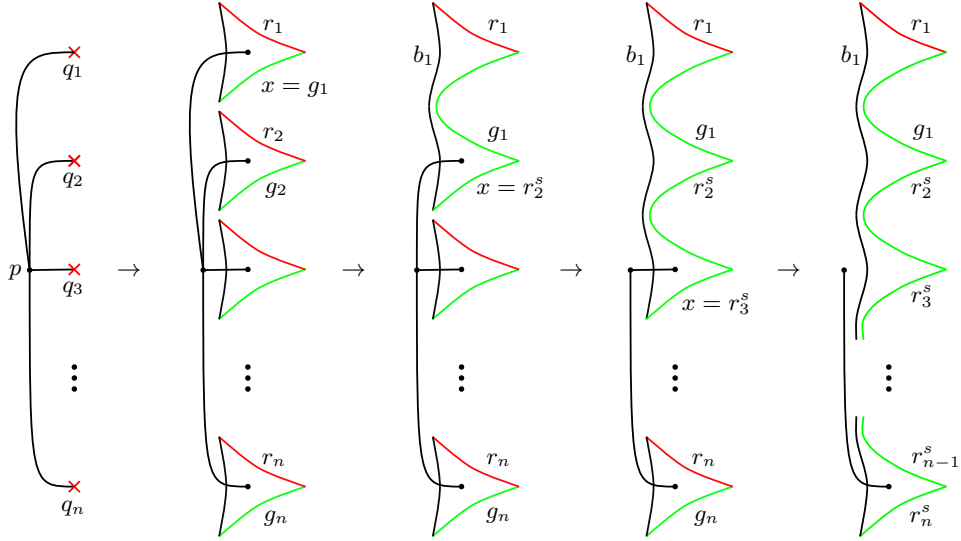
\begin{figure}\centering
	\labellist
	\small\hair 2pt
	\pinlabel $p$ at 0 110
	\pinlabel $q_1$ at 24 193
	\pinlabel $q_2$ at 24 148
	\pinlabel $q_3$ at 24 103
	\pinlabel $q_n$ at 24 12
	\pinlabel $\rightarrow$ at 47 110
	\pinlabel $r_1$ at 107 212
	\pinlabel $x=g_1$ at 116 186
	\pinlabel $r_2$ at 107 167
	\pinlabel $g_2$ at 108 144
	\pinlabel $r_n$ at 107 32
	\pinlabel $g_n$ at 108 8
	\pinlabel $\rightarrow$ at 140 110
	\pinlabel $r_1$ at 200 212
	\pinlabel $b_1$ at 169 200
	\pinlabel $g_1$ at 200 167
	\pinlabel $x=r_2^s$ at 205 145
	\pinlabel $r_n$ at 200 32
	\pinlabel $g_n$ at 200 9
	\pinlabel $\rightarrow$ at 230 110
	\pinlabel $r_1$ at 285 212
	\pinlabel $b_1$ at 257 200
	\pinlabel $g_1$ at 286 168
	\pinlabel $r_2^s$ at 285 145
	\pinlabel $x=r_3^s$ at 290 97
	\pinlabel $r_n$ at 285 32
	\pinlabel $g_n$ at 286 8
	\pinlabel $\rightarrow$ at 320 110
	\pinlabel $r_1$ at 375 212
	\pinlabel $b_1$ at 346 200
	\pinlabel $g_1$ at 376 168
	\pinlabel $r_2^s$ at 375 145
	\pinlabel $r_3^s$ at 375 100
	\pinlabel $r_{n-1}^s$ at 381 33
	\pinlabel $r_n^s$ at 375 9
	\endlabellist
	\includegraphics[width=.75\linewidth]{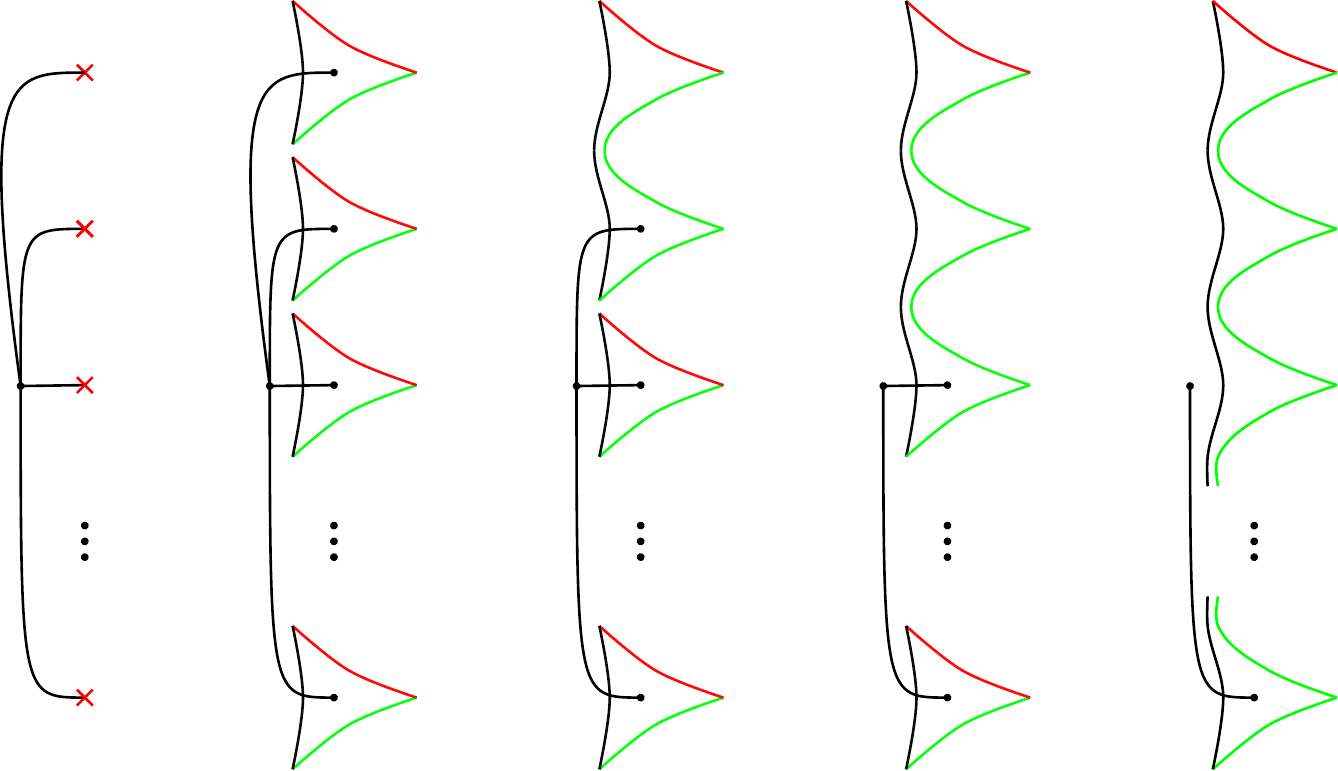}
	\caption{Wrinkling the Lefschetz critical points and then merging the resulting triangles results in a crown map, whose vanishing cycles form a coronation of the Lefschetz fibration. The main circle is always the top circle here.}
	\label{mergesfig}
\end{figure} Once $i-1$ triangles have been merged into the main circle, the main circle will have $i+1$ vanishing cycles, and the method of proof is to show that merging the $i^{th}$ triangle into the main circle as below results in the main circle having a sequence of vanishing cycles which is by definition a coronation of $l_1,\ldots,l_i$. This sequence is then shown to be slide-equivalent to a pseudocoronation of $l_1,\ldots,l_i$. Thus, inductively, there is a coronation which is slide-equivalent to a pseudocoronation of $l_1,\ldots,l_n$. Then, according to \cite[Theorem 4.9]{W3}, all genus-$g$ crown maps which are homotopic to the original Lefschetz fibration have crown diagrams which are related by the genus-preserving moves \emph{handleslide, multislide} and \emph{shift}, which by \cite[Theorem 4.3]{W3} and Theorem \ref{thm1} preserve the slide-equivalence class of the diagrams. For this reason, \emph{all} possible coronations are slide-equivalent to some pseudocoronation. 

To prove that all pseudocoronations of a given sequence are slide-equivalent, note that one can form any pseudocoronation by choosing a smooth embedded path $\phi$ starting at a point of $r_1$ and transversely intersecting $r_1,g_1,\ldots,r_{n-1}$ in order, with no other intersections, and ending at $r_n$, then replace a neighborhood of this path with Figure \ref{wrinklefig3}. In that figure, it is a short exercise to see that one may slide $r_1,\ldots,r_n$ so that they may be isotoped off of the pictured region, then one may change the choice of $\phi$ and slide them back on.

It remains to prove that merging the $i^{th}$ triangle into the main circle results in the desired slide-equivalence class. The first step is to derive the resulting sequence of vanishing cycles. Clearly the first triangle's surgered diagram is that of a coronation of the single Lefschetz vanishing cycle $l_1$, and it is slide-equivalent to itself. Now assume the first $i-1$ triangles have been merged into the main circle for some $i>1$, and that the sequence of vanishing cycles of the surgered diagram for the main circle is slide-equivalent to the coronation \begin{equation}\label{coro1}b_1,r_1,g_1,r_2,r_3,\ldots,r_{i-1}.\end{equation} 
To merge the $i^{th}$ triangle into the main circle, begin with the surgered diagram for the $i^{th}$ triangle. Decorate it further with the disks corresponding to $b_1^1,b_1^2$ and the last vanishing cycle $x$ of the main circle. By the induction hypothesis, $x$ is slide-equivalent to the corresponding element of a pseudocoronation as an element of the surgered diagram for the main circle, so denote it $g_1$ if $i=2$, or by $r_{i-1}^s$ if $i>2$ (See Figure \ref{mergexfig1} for an example).
\begin{figure}
	\begin{subfigure}[t]{0.30\textwidth}
		\centering
		\labellist
		\small\hair 2pt
		\pinlabel $r_i$ at 99 16
		\pinlabel $x$ at 125 16
		\pinlabel $b_1^2$ at 124 34
		\pinlabel $b_1^1$ at 124 65
		\pinlabel $b_i^2$ at 95 34
		\pinlabel $b_i^1$ at 95 65
		\pinlabel $g_i$ at 82 55
		\endlabellist
		\includegraphics[width=\textwidth]{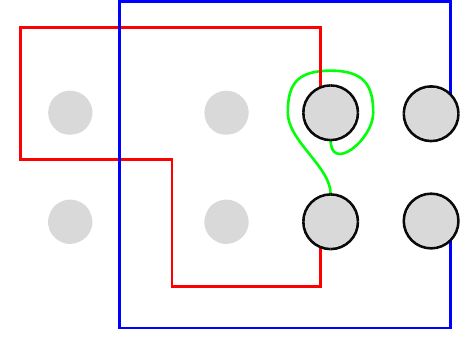}
		\caption{Initial diagram in which $x$ is the last vanishing cycle of the main circle. If $i=2$, then $x=g_1$. If $i>2$, then $x=r_{i-1}^s$.}
		\label{mergexfig1}
	\end{subfigure}
	\hspace{0.04\textwidth}
	\begin{subfigure}[t]{0.30\textwidth}
		\centering
		\labellist
		\small\hair 2pt
		\pinlabel $r_i^s$ at 94 16
		\pinlabel 2 at 90 34
		\pinlabel 1 at 90 65
		\endlabellist
		\includegraphics[width=\linewidth]{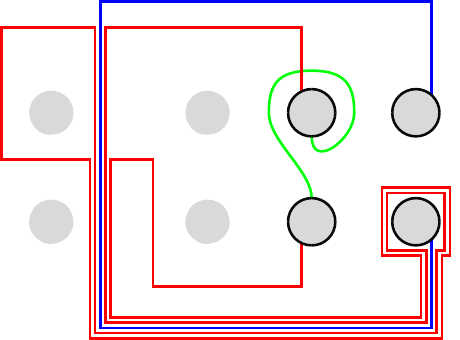}
		\caption{Slide $r_i$ over $b_1$. Write superscripts 1,2 for $b_i^1$, $b_i^2$.}
		\label{mergexfig2}
	\end{subfigure}
	\hspace{0.04\textwidth}
	\begin{subfigure}[t]{0.30\textwidth}
		\centering
		\labellist
		\small\hair 2pt
		\pinlabel 1 at 90 34
		\pinlabel 2 at 90 65
		\pinlabel $r_i^s$ at 94 16
		\pinlabel $g_i^s$ at 81 52
		\endlabellist
		\includegraphics[width=\linewidth]{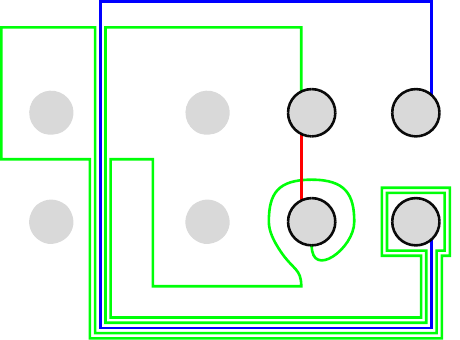}
		\caption{Switching push map: Push $b_i^1$ along $r_i^s$ to $b_i^2$, push $b_i^2$ up to where $b_i^1$ was. It turns out that $g_i$ has become slide-equivalent to $r_i$ and vice-versa, so labels are switched.
		}
		\label{mergexfig3}
	\end{subfigure}\vspace{.5cm}
	\begin{subfigure}[t]{0.30\textwidth}
		\centering
		\labellist
		\small\hair 2pt
		\pinlabel 2 at 99 31
		\pinlabel 1 at 90 65
		\endlabellist
		\includegraphics[width=\linewidth]{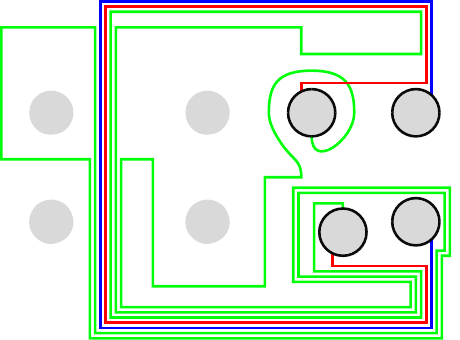}
		\caption{Switching push map: Push $b_i^2$ along $x$ and push $b_i^2$ up to where $b_i^1$ was.}
		\label{mergexfig4}
	\end{subfigure}\hspace{0.1\textwidth}
	\begin{subfigure}[t]{0.30\textwidth}
		\centering
		\labellist
		\small\hair 2pt
		\pinlabel $r_i^s$ at 138 65
		\endlabellist
		\includegraphics[width=\linewidth]{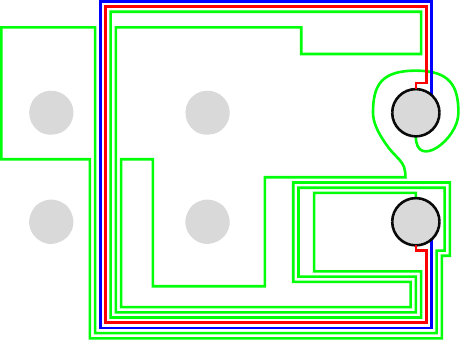}
		\caption{Move $b_i^j$ to coincide with $b_1^j$.}
		\label{mergexfig5}
	\end{subfigure}
	\caption{An example of merging a triangle into the main circle using the technique from Section \ref{cuspmergesection}. These modifications are a particular isotopy of the surgered diagram for the $i^{th}$ triangle which causes $r_i$ and $b_i^j$ to coincide with $x$ and $b_1^j$, respectively. The remaining vanishing cycle $g_i$ becomes $r_i^s$, the last vanishing cycle associated to the main circle. Disks without boundary circles are included as indication of the genus of the surface, and all are paired vertically.}\label{mergexfig}
\end{figure} 

Following the method of Section \ref{cuspmergesection}, the next step is to isotope the surgered diagram for the $i^{th}$ triangle until its disks coincide with those of the main circle, and until $r_i$ coincides with $x$. This is done in the following steps: First, arbitrarily orient $x$ so that $r_i\cap x$ becomes an ordered set. As in Figure \ref{mergexfig2}, use the part of $x$ that lies after the last element of $r_i\cap x$ as a path along which $r_i$ slides over $b_1$, thereby removing the last element of $r_i\cap x$. Repeat until $r_i\cap x$ is empty. 

Continuing the isotopy, the next step is to apply what this paper will call a \emph{switching push map} to the diagram for the $i^{th}$ triangle. A switching push map is a point pushing map in the sense of the Birman exact sequence, except that the map is realized as an isotopy of the embedding of each of the two disks being pushed instead of a single point, and the two disks switch places: One disk follows a vanishing cycle connecting the two disks, while the other disk travels along another such vanishing cycle to the starting point of the first disk. One last aside before describing the switching push maps: Note that Figure \ref{mergexfig} has the two pairs of disks nearby to each other as opposed to on top of each other as in Figure \ref{wrinklefig3} purely for the sake of distinguishing which disks are moving around, and the following switching push maps should be better understood with the pairs of disks superimposed like in Figure \ref{wrinklefig3}, so that each push map is a genuine switching push map (consider how, for example, the change from Figure \ref{mergexfig3}-\ref{mergexfig4} would only be a genuine push map if the pairs of disks coincided perfectly). Now, for the first switching push map, push the disk bounded by $b_i^2$ along $r_i$ and push $b_i^1$ along the short vertical linear path to where $b_i^2$ started, causing $r_i^s$ to become a short linear path as in Figure \ref{mergexfig3}. Note this short linear vertical path is slide-equivalent to $g_1$. Then apply another switching push map as in Figure \ref{mergexfig4}, in which the disk bounded by $b_i^2$ travels along $x$, so that $r_i^s$ becomes parallel to $x$. Finally, slightly perturb the diagram for the last triangle so that its disks coincide exactly with those of the main circle, and $r_i^s$ coincides with $x$ as in Figure \ref{mergexfig5}. At this point, $b_i^1$, $b_i^2$, the disks they bound, and the parallel copy of $x$ are deleted from the complete surgered diagram of the whole Morse 2-function, and the vanishing cycle $r_i^s$ in Figure \ref{mergexfig5} becomes the last vanishing cycle for the main circle. 

\begin{figure}
\begin{subfigure}[t]{0.20\textwidth}
	\centering
	\labellist
	\small\hair 2pt
	\pinlabel $a$ at 35 13
	\pinlabel $b$ at 4 5
	\pinlabel $c$ at 71 45
	\endlabellist
	\includegraphics[width=\textwidth]{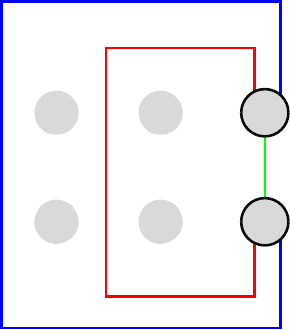}
	\caption{Beginning configuration.}
	\label{slidecheckfig1}
\end{subfigure}\hspace{0.04\textwidth}
\begin{subfigure}[t]{0.20\textwidth}
	\centering
	\includegraphics[width=\linewidth]{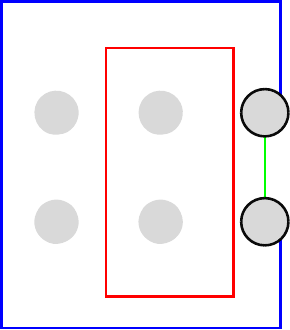}
	\caption{Slide $a$ over $c$.}
	\label{slidecheckfig2}
\end{subfigure}\hspace{0.04\textwidth}
\begin{subfigure}[t]{0.20\textwidth}
	\centering
	\labellist
	\small\hair 2pt
	\pinlabel $a^s$ at 37 15
	\endlabellist
	\includegraphics[width=\linewidth]{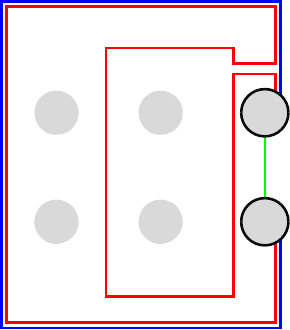}
	\caption{Slide $a$ over $b$.}
	\label{slidecheckfig3}
\end{subfigure}\hspace{0.04\textwidth}
\begin{subfigure}[t]{0.20\textwidth}
	\centering
	\labellist
	\small\hair 2pt
	\pinlabel $a^s$ at 37 15
	\endlabellist
	\includegraphics[width=\linewidth]{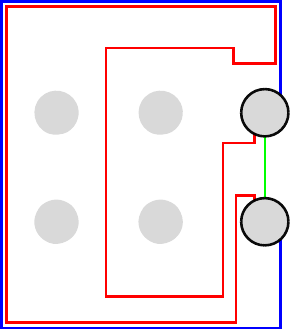}
	\caption{Output of the switching push map.}
	\label{slidecheckfig4}
\end{subfigure}
\caption{Applying a switching push map to a vanishing cycle $a$ according to $b$ and $c$ results in a vanishing cycle $a^s$ which may also be obtained by sliding $a$ over $b$ and then over $c$. Figures \ref{slidecheckfig3} and \ref{slidecheckfig4} differ by an isotopy of $a^s$.}\label{slidecheckfig}
\end{figure}

Having described the calculation of the vanishing cycle $r_i^s$ that is added to that of the main circle when the $i^{th}$ triangle is merged into the main circle, the goal is to prove that $b_1,r_1,g_1,r_2,r_3,\ldots,r_i^s$ is indeed slide-equivalent to $b_1,r_1,g_1,r_2,r_3,\ldots,r_i$. Figure \ref{slidecheckfig1} depicts the general starting position, where the two disks are connected by three pairwise disjoint vanishing cycles $a,b,c$. In this figure, $c$ represents the result of sliding $g_1$ over $b_1$, so sliding over $c$ is therefore tantamount to sliding over $g$ and then sliding over $b_1$, so sliding $r_i$ over such a vanishing cycle preserves the slide-equivalence class of $b_1,r_1,g_1,r_2,r_3,\ldots,r_i$. Since $b$ represents the vanishing cycle $r_i^s$ in Figure \ref{mergexfig2}, which is obtained by repeatedly sliding $r_i$ over $b_1$, sliding over such a vanishing cycle for the first switching push map similarly preserves the slide-equivalence class of $b_1,r_1,g_1,r_2,r_3,\ldots,r_i$. Finally, for the second switching push map, $b$ corresponds to $x$, with similar results. Overall, the sequence of vanishing cycles over which $g_i$ must slide to become the last vanishing cycle for the main circle reads $g_1^s$, then $r_i^s$ (at which point it is slide-equivalent to $r_i$), then $g_1^s$ and finally $x$.
\end{proof}

\bibliographystyle{plain}

\end{document}